\documentclass[10pt]{amsart}
\usepackage{amsmath,amssymb,amscd,enumerate,bbm}
\usepackage{cite}

\renewcommand{\MR}[1]{}

\newcommand{\titl}{THE CANONICAL STRIP, I}

\title{{\titl}}
\author{V. Golyshev}
\date{}

\def\R{{\Bbb{R}}}

\def\OO{{\mathcal{O}}}
\renewcommand{\phi}{{\varphi}}

\newcommand{\dimens}{{\mathrm{dim \;} }}

\newcommand{\Real}{\mathop{\mathrm{Re}}}

\newcounter{pphcounter}[section]
\renewcommand{\thepphcounter}{\thesection.\arabic{pphcounter}}
\newcommand{\pph}{\bigskip \refstepcounter{pphcounter}
    \bf  \thepphcounter. \rm}

\makeatletter
\@addtoreset{equation}{pphcounter}
\makeatother

\renewcommand{\proof}{\bigskip \bf Proof. \rm}

\def\A1{{{\Bbb{A}}^1}}
\def\P1{{{\Bbb{P}}^1}}

\def\SL2{{\mathrm SL2}}

\renewcommand{\P}{{\Bbb{P}}}

\begin{document}

\begin{center}
\maketitle


\bigskip

\end{center}

\bigskip

\parbox{340pt}{\small \bf Abstract. \rm We introduce a canonical
strip hypothesis for Fano varieties. We show that the canonical strip
hypothesis for a Fano variety implies that the zeros of the Hilbert
polynomial of embedded Calabi--Yau
and general type hypersurfaces are located on a vertical line.
This extends, in particular, Villegas's `polynomial
RH' for intersections in projective spaces to the
case of CY and general type hyperplane sections in Grassmannians.
We state a few conjectures on the Ehrhart polynomials 
of certain fan polytopes.}

\bigskip
\bigskip
\bigskip

\section{Vanishing theorems and the canonical strip}

The prototypal vanishing
theorem is due to Kodaira. Let $X$ be a smooth complex projective
variety of dimension $n$, and let $K_X$ be its canonical divisor.
Then higher cohomology groups
$H^i(K_X+L), \: i>0$, vanish for any ample line bundle $L$.
It is usually viewed as a theorem in Kahler geometry; in fact, the Akizuki--Nakano
theorem says that for a positive hermitian holomorphic line bundle
$E$ on a compact Kahler
manifold $X$ one has $H^{p,q}(X,E)=0$ for $p+q \ge n+1$; the previous assertion
follows if one sets $p=n$. A few semipositive and partially positive versions
are discussed in Demailly's textbook \cite{Demailly07}. As a generalisation for vector bundles,
one has Nakano's vanishing theorem: if a hermitian vector bundle $E$ is
\emph{Nakano positive},
then $H^{n,q}(X,E)=0$ for $q>0$. Back to algebraic geometry and line
bundles, strong and subtle generalisations
of the Kodaira vanishing theorem have been proved, cf. \cite{EV92}.
The import of Kodaira's theorem varies according to which world, Fano
or Calabi--Yau, or general type, we are in. We are mainly interested
in Fano in this note; for a Fano variety
(i.e. in the case when $K_X<0$), the corollary is that the
line bundles in the
\emph{canonical strip}, $K_X < E < 0$, are acyclic.

Another Fano--pertaining subject which involves the
consideration of the canonical strip,
is theory of helices and stability conditions. For a Fano $X$ with a full
\emph{exceptional collection} $\langle E_i\: \mathrm{'s} \rangle$ that
has $\OO$ for its member,
a turn of the respective helix between
$-K_X$  and $\mathcal O$ (non--inclusive) consists of acyclic objects.
The issue of Nakano positivity (or Griffiths positivity) of the objects
in such turn twisted by $-K_X$ has not been addressed, to our knowledge.

The qualitative phenomena referred
to above, the acyclicity of the line bundles in the canonical strip of a Fano variety,
and the acyclicity
of the vector bundles in a canonical turn, appear to admit quantitative treatment via
the study of location of zeros
of the Hilbert polynomial. Our theorem \ref{theo-RH}
was inspired by the elegant note \cite{Rodriguez-Villegas02} by Rodriguez--Villegas.

\pph \bf The canonical strip and the canonical line hypotheses. \rm
Let $H(z)$ be the Hilbert polynomial
of a [Fano or general type]
variety $X$, so that $H(n)= H_{-K_X}(n)= \chi (n (-K))$
for integral $n$'s.

(CS)We say that $X$ satisfies \emph{the canonical strip hypothesis}
if all roots $z_i$ of
$H(z)$ are in the \emph{canonical strip} $-1 < \Real z <0$.

(NCS)We say that $X$ satisfies \emph{the narrowed canonical strip hypothesis}
if all roots $z_i$ of
$H(z)$ are in the \emph{narrowed canonical strip}
$$-1+\frac{1}{\dimens X +1} \le \Real z \le -\frac{1}{\dimens X +1}.$$

(CL) We say that $X$ satisfies \emph{the canonical line hypothesis},
if all roots $z_i$ are on the vertical line
$\Real z = -1/2 $.

It is clear that (CL) $\Longrightarrow$ (NCS) $\Longrightarrow$ (CS).

\pph \bf Calabi--Yau: an embedded version. \rm
For a Calabi--Yau type $X$ the Hilbert polynomial
with respect to the anticanonical class is not too exciting. We consider instead
embedded Calabi--Yaus. We say that a Calabi--Yau
$X$ embedded as an anticanonical section in a Fano $F$ satisfies
the canonical line hypothesis
if the roots of its Hilbert polynomial with respect to the restriction
of $-K_F$ to $X$ are purely imaginary.

\pph \bf Curves. \rm
For a genus $g$ curve, $\chi (z (-K))=(2-2g)(z+1/2)$,
and the canonical line hypothesis holds. For an elliptic curve, embedded as a cubic
in $\P^2$, one has $\chi (\OO (9z))=9z$, and the canonical line hypothesis holds.

\pph \bf Surfaces. \rm The Riemann--Roch--Hirzebruch formula yields, for surfaces,
$$H(z)=1/2\,{z}^{2}{{ c_1}}^{2}+1/2\,{{ c_1}}^{2}z+1/12\,{{ c_1}}^{2}+1
/12\,{ c_2}.$$

The two roots
of $H(z)$ are
$$-1/2 \pm 1/6\sqrt{\frac{-6 c_2+3c_1^{2}}{c_1^2}}.$$

Thus, one has:
$$
\begin{array}{lcl}
X  \text{ satisfies (CS)} & \Longleftrightarrow  & c_1^2 \ge -c_2 \\
X  \text{ satisfies (NCS)} & \Longleftrightarrow & c_1^2 \le 3c_2 \\
X  \text{ satisfies (CL)} & \Longleftrightarrow & c_1^2 \le 2c_2
\end{array}
$$
For del Pezzo surfaces, the maximal possible value of
$\displaystyle{{3- \frac{6 c_2}{c_1^2}}}$
is $1$; it is attained on $\P ^2$. Thus, (NCS) holds for del Pezzos. For
surfaces of general type, (NCS) holds by Yau.


For an embedded K3's polarized by a class $h$,
$\chi (\OO (h))=1/2h^2z^2+2$, and the canonical line hypothesis holds.

\pph \bf Threefolds. \rm  Riemann--Roch--Hirzebruch says now that
$$H(z)=1/6\,{z}^{3}{{  c_1}}^{3}+
1/4\,{{  c_1}}^{3}{z}^{2}+\left (1/12\,{{  c_1}}^{2}+
1/12\,{  c_2}\right )z{  c_1}+1/24\,{  c_1}{  c_2},$$

and the three roots are
$$-1/2, -1/2 \pm \frac{1}{2} \sqrt{\frac{c_1^3-2c_1c_2}{c_1^3}}.$$

One has:
$$
\begin{array}{lcl}
$$X  \text{ satisfies (CS)} & \Longleftrightarrow & \dfrac{-2c_1c_2}{c_1^3} < 0 \\
X  \text{ satisfies (NCS)} & \Longleftrightarrow & \dfrac{-2c_1c_2}{c_1^3} \le -3/4\\
X  \text{ satisfies (CL)} & \Longleftrightarrow  & \dfrac{-2c_1c_2}{c_1^3} \le -1
\end{array}
$$

In the Fano case, the vanishing theorem gives $c_1c_2=24$, so
the maximal possible value of $-{\dfrac{2c_1c_2}{c_1^3}}$
is $-3/4$, and it
is attained on $\P ^3$, so (NCS) is true.
The formula for the roots  shows that the narrowed canonical
strip hypothesis holds for minimal threefolds of general type: Yau's result
for threefolds is $c_1^3 \ge 8/3 c_1c_2$.

\pph \bf Grassmannians. \rm
Projective spaces satisfy the narrow canonical strip hypothesis, the zeros
being $\dfrac{-i}{n+1}, i=1,\dots , n$. This generalizes to Grassmannians,
as shown by Hirzebruch. Let $G(k,N)$ be the Grassmannian of $k$--planes
in an $N$--dimensional space. Assume that $N \ge 2k$. Let $\phi (x)$ be the piecewise
linear function of real argument $x$ given by
$\phi (x)= \min (k,-x,x+n+1).$
Then  \cite{Hirzebruch58}
$$H(z)=c \prod_{i=1}^{n} (z+\frac{i}{n+1})^{\phi(-i)}.$$

\section{The canonical line hypothesis for embedded varieties}

As we saw above, varieties of general type need not satisfy (CL).
The situation is different with embedded varieties of general type.

\pph \label{theo-RH} \bf Theorem. \rm Let $F$ be a Fano variety which satisfies the
canonical strip hypothesis, and let $X$ be its  general type (resp. Calabi--Yau
type) section in the linear system $-nK_X$, $n>1$ (resp. $n=1$).
Then (CL) holds for the embedded variety $X$.

\proof Let $H_F(z)$ be the Hilbert polynomial of $F$,
and let $H_r(z)$ denote the Hilbert polynomial of $X$ with respect to the restricted
$-K_F$. For general type $X$, the adjunction formula says that the anticanonical
class of $X$ is a multiple of the restricted anticanonical class of $F$, so we
may prove the vertical line statement for $H_r$ instead.
Then
$$H_r(z)=H_F(z)-H_F(z-n).$$
The fact that the zeros of $H_r(z)$  are on the line $\Real z =\dfrac{n-1}{2}$
follows from the

\pph \bf Lemma. \rm  Let $H(z) \in \R [z]$ satisfy

i) $H(-1-z)=\pm H(z)$,

ii) all roots  $z_i$ of $H$ are [strictly] in the left half--plane.

Then, for any real $s \ge 1$ and for all roots $\zeta_j$  of $H(z)-H(z-s)$,
one has $\Real z=\dfrac{s-1}{2}$.

\proof Assume there is a root $\zeta=\zeta_j$ of $H_r (z)$
with $\Real \zeta < \dfrac{s-1}{2}$.
By assumption,
$$\prod (\zeta-z_i)=  \prod(\zeta-z_i-s). $$
Let $\mu_s$ be the map that takes $z$ to $(s-1)-\bar z.$ It particular,
it establishes a 1-1 mapping between the factors in the products above:
$$\prod (\zeta-z_i)=  \prod(\zeta-\mu_s(z_i)). $$\
However, for any given $i$
$$\mid \zeta-z_i \mid < \mid \zeta-\mu_s(z_i)\mid,$$
as
$\zeta$ is to the left of the axis of the symmetry $\mu$.

The case $\Real \zeta > \dfrac{s-1}{2}$ is treated in the same manner.
This contradiction proves Lemma and Theorem \ref{theo-RH}.

\pph \bf Corollary. \rm A section of $-mK_X, \: m>0$ of a Grassmann variety
$X$ satisfies the canonical line hypothesis.
\qed

\pph \bf Problems. \rm

 \bf A. \rm  Are there Fano or general type varieties
that do not satisfy the canonical strip hypothesis? How can one characterize
those varieties that satisfy (CS) but not (NCS)?

\bf B. \rm  Study zeros of the Ehrhart polynomials \cite{BDD+05}, \cite{BHW07} of
Fano polytopes as compared to the Ehrhart polynomials of non--Fano polytopes.
Is there a  characterization of the former? We conjecture that the
Ehrhart polynomial of a fan polytope in the space of cocharacters
of a smooth toric Fano in dimensions $1,2,3,4,5$ has all zeros in the line $\mathop{\mathrm{Re}} z =-1/2$ 
\footnote{Dimensions $1$ and $2$ are easy. 
C. Shramov informed me recently that he had proved this
in dimension $3$.}. We furthermore conjecture that the Ehrhart polynomial of a fan polytope of any terminal Gorenstein toric Fano 3--fold (which, by \cite{Friedman86}, admits a smoothing)
has the same property.

\bf C. \rm Starting with dimension $4$, the classical Routh--Hurwitz
stability criterion furnishes a set of polynomials in Chern numbers
whose positivity/non--negativity implies (CS) and (NCS).
Is it possible to prove inequalities of such type using Yau's theorem,
or a hierarchy of enhancements is needed?

\bigskip
\bigskip

\bf Acknowledgements. \rm I thank Don Zagier for referring me to Villegas's work, and Tyrrell McAllister for pointing me to \cite{BHW07}.
I thank Friedrich Hirzebruch for the discussion of the subject, and Ivan Cheltsov
 and Constantin Shramov for their advice.

\nocite{GS07}

\bigskip
\bigskip
\bigskip
\bigskip
\bigskip

\bibliographystyle{amsalpha}

\end{document}